\documentclass[twoside,12pt]{article}
\makeatletter
\newcommand{\comment}[1]{}
\def\ifundefined#1{\expandafter\ifx\csname#1\endcsname\relax}
\makeatletter
\usepackage{theorem}
     \theorembodyfont{\slshape}
        \newtheorem{thm}{Theorem}[section]
     
     \newtheorem{lem}[thm]{Lemma}

     \theorembodyfont{\upshape}
     \newtheorem{defn}[thm]{Definition}
     
     \newtheorem{example}[thm]{Example}

     \newtheorem{rem}[thm]{\mdseries\scshape Remark}
\newenvironment{proof}[1][\proofname]{\par
  \normalfont
  \topsep6\p@\@plus6\p@ \trivlist
  \item[\hskip\labelsep\scshape
    #1.]\ignorespaces
}{%
  $\qed$\endtrivlist
}
\newcommand{\proofname}{Proof}

\newcommand{\epigraph}[3]{\par
\hfill\parbox{0.4\textwidth}{\footnotesize #1 \par \hfil #2 \it
#3}\par\medskip}
\providecommand{\dedicatory}[1]{}
\providecommand{\keywords}[1]{\begingroup \def \protect {\noexpand \protect
\noexpand }\xdef \@thefnmark { }\endgroup \@footnotetext{{\em Keywords and
phrases.\/} #1}}
\providecommand{\AMSMSC}[2]{\begingroup \def \protect {\noexpand \protect
\noexpand }\xdef \@thefnmark { }\endgroup \@footnotetext{{1991 \it
Mathematical Subject Classification.\/} Primary: #1; Secondary: #2.}}

\newcommand{\authorshort}[1]{\gdef\@authorshort{#1}}
\newcommand{\titleshort}[1]{\gdef\@titleshort{#1}}
\gdef\oldmaketitle{\par \begingroup \renewcommand \thefootnote
{\fnsymbol {footnote}}\def \@makefnmark {\hbox to\z@
{$\m@th ^{\@thefnmark }$\hss }}\long \def \@makefntext ##1{\parindent
1em\noindent
\hbox to1.8em{\hss $\m@th ^{\@thefnmark }$}##1}\if@twocolumn \ifnum
\col@number =\@ne \@maketitle \else \twocolumn [\@maketitle ]\fi
\else \newpage \global \@topnum \z@ \@maketitle \fi
\thispagestyle {plain}\@thanks \endgroup \setcounter {footnote}{0}
\let\thanks \comment \relax \gdef \@thanks {}}

\def\@cite#1#2{{#1\if@tempswa , #2\fi }}
\def\@citex[#1]#2{\let \@citea \@empty \@cite {\@for \@citeb :=#2\do
{\@citea \def\@citea {,\penalty \@m \ %
}\edef \@citeb {\expandafter \@iden \@citeb
}\if@filesw \immediate\write\@auxout{\string \citation {\@citeb }}\fi
\@ifundefined
{b@\@citeb }{{\reset@font \bfseries ?}\G@refundefinedtrue
\@latex@warning
{Citation `\@citeb ' on page \thepage \space undefined}}%
{\hbox {[\csname b@\@citeb \endcsname]}}}}{#1}}
\makeatletter
\newcommand{\reportenum}[3][]{\gdef\rep@rtenum{#2}
\gdef\rep@rteyear{#3}\gdef\wh@reappear{#1}}
\let\@ldmaketitle=\maketitle
\renewcommand{\maketitle}{{\def\newpage{}
{\scriptsize\parbox[t]{0.3\textwidth}{\noindent Reporte Interno
\# \rep@rtenum\\
Departamento de Matem\'aticas\\CINVESTAV del IPN\\Mexico City, \rep@rteyear}
\hfill
\parbox[t]{0.5\textwidth}{\wh@reappear}}
\@ldmaketitle}}
\makeatother
     \theorembodyfont{\upshape}
     \newtheorem{procedure}[thm]{Procedure}

\newcommand{\spt}{\object{supp}}

\newcommand{\smnegsk}{\mkern-3mu}
\usepackage{amsfonts}
\reportenum[To appear in\\{\em Mathematical Methods in\\Applied
Sciences}]{166}{1994}
\newcommand{\algebra}[1]{\ensuremath{{\mathfrak #1}}}

\newcommand{\Cliff}[2]{\ensuremath{{\bf Cl}(#1,#2)}}
\newcommand{\object}[1]{\ensuremath{\mathrm{#1}}\,}
\newcommand{\Space}[2]{\ensuremath{ {{\mathbb #1}^{#2}} }}

\newcommand{\FSpace}[2]{{\ensuremath{ #1_{#2} }}}

\newcommand{\such}{\ensuremath{\,|\,}}
\newcommand{\bos}{{\mathcal B}}
\ifundefined{qed}
    \DeclareMathSymbol{\qed}{0}{AMSa}{"03}
\fi
\newcommand{\Cstar}{$C^{*}$}

\newcommand{\norm}[1]{\parallel\smnegsk#1\smnegsk\parallel}
\newcommand{\modulus}[1]{\mid\smnegsk#1\smnegsk\mid}
\newcommand{\scalar}[2]{\langle #1,#2\rangle}

\providecommand{\eqref}[1]{\textup{(\ref{#1})}}

\authorshort{Vladimir V. Kisil and Enrique Ram\'{\i}rez de Arellano}
\titleshort{The Riesz-Clifford Functional Calculus}

\begin{document}
\title{The Riesz-Clifford Functional Calculus \\
for Non-Commuting Operators \\
and Quantum Field Theory\thanks{This work was partially
supported
by CONACYT
Project 1821-E9211, Mexico.}}

\author{Vladimir V. Kisil\thanks{On leave from the Odessa State
University.}\\
\normalsize  e-mail:\ttfamily vkisil@mvax1.red.cinvestav.mx \and
Enrique Ram\'{\i}rez de Arellano\\
\normalsize  e-mail:\ttfamily eramirez@mvax1.red.cinvestav.mx \and
Departamento de Matem\'aticas,
CINVESTAV del I.P.N.,\\
Apartado Postal 14-740,
07000, M\'exico, D.F. M\'exico  \\
fax: (525)-752-64-12}
\date{November 24, 1994}
\maketitle

\begin{abstract}
We present a Riesz-like hyperholomorphic functional calculus for a
set
of non-commuting
operators based on the Clifford analysis. Applications to the quantum
field theory are described.
\keywords{Functional calculus, Weyl calculus, Riesz calculus, Clifford
analysis, quantization, quantum
field theory.}
\AMSMSC{47A60}{81T10}
\end{abstract}
\newpage
\tableofcontents
\section{Introduction}
\epigraph{The quantum mechanical connection was certainly the
strongest of the motivations for the development of operator
algebra.}{Irving E. Segal,~\cite{Segal94}}{}
The quantization
procedure~\cite{Berezin74,Coburn94a,Folland89,Kirillov90} requires
the construction of some mapping from an algebra of real valued functions
(the algebra of classic observables) to an algebra of self-adjoint
operators (the algebra of quantum observables). Analogously, one
should be
able (in some sense) to construct functions from self-adjoint
(non-commuting) operators. The importance of this problem in the
 physics and the non-triviality of the corresponding
mathematical problems has been amply illustrated during the last 50 years (in
addition to the above mentioned papers see also recent
ones~\cite{Effros94,Rieffel94,Segal94} and their rich bibliographies).

The problem of quantization may be reformulated mathematically as a
problem on construction of a functional calculus. The already classical
holomorphic calculus, which is based on analysis of several
complex
variables,
allows only the situation of mutually commuting operators (see the
milestone
papers of J.~L.~Taylor~\cite{JTaylor70,JTaylor72}). Other approaches
also lead only to the commuting case (see, for
example,~\cite{McInPryde87}).
Such
constructions do not meet the needs of quantum mechanics, where
non-commutativity of observables is a principal assumption.

The main
goal of
the present paper is a construction of a Clifford hyperholomorphic calculus of
{\em several
non-commuting\/} operators and to show its applicability to
quantum field theory.

It should be mentioned that functional calculus may be considered at
two different levels. The first one is the original problem of quantum
mechanics of constructing functions from the operators of coordinates
$Q_j$
and impulses $P_j$, which satisfy the well known Heisenberg commutation
relations
\begin{equation}\label{eq:heisen-comm}
[Q_j,P_i]=i\hbar\delta _{ij}I.
\end{equation}
In this context there are
many different approaches to quantization with a large number of nice
results (for example, the PDO
calculus~\cite{Folland89,Hormander85,Shubin87,MTaylor81} and the
Toeplitz
operators in the Segal-Bargman space~\cite{Coburn94a,CobXia94}). The
richness of
the calculus is rooted in the structure of the Heisenberg
group~\cite{Howe80a,Howe80b}. At a second level  one can construct a
functional calculus from an arbitrary finite set of self-adjoint
operators $\{T_j\}, 1\leq j\leq m$. In this case the number of
possible
quantizations and the precise results obtained is much smaller~\cite{Maslov73}
and
sometimes only a more-or-less formal power series may be established
in such a setting.
In~\cite{Kisil94e} it was shown that one can preserve the essential
features at the first level if the assumption is made, that the operators
$\{T_j\},
1\leq j\leq m$ represent some Lie algebra \algebra{g}.

 In Section~\ref{se:approach} of this work we
discuss
general properties and known results of functional calculus.
In
Section~\ref{se:riesz} we construct the Riesz-like  hyperholomorphic
functional calculus for a set of non-commuting operators based on
Clifford analysis. To extend the calculus from commuting operators to
non-commuting ones we restrict the family of functions from smooth ones to
hyperholomorphic functions. It should be noted that hyperholomorphic
functions do not form an algebra (at least under usual multiplication) and this
is reflected at Definition~\ref{de:calculus} of Riesz-Clifford calculus.
 In Section~\ref{se:field} it is shown that the
Riesz-Clifford functional calculus is a relevant model in quantum
field
theory corresponding to the Weyl calculus in quantum mechanics.

We consider only the case of self-adjoint operators on the Hilbert
space
$H$ (or of elements in an abstract \Cstar-algebra). It is also possible
to
consider the more general setting of self-adjoint operators in a Banach
space~\cite{Anderson69} (or of elements in an abstract Banach algebra).

We are grateful to V.~V.~Kravchenko, J.~Ryan, and
I.~Spit\-kov\-sky for their valuable comments.

\section{Different Approaches to Functional
Calculus}\label{se:approach}
In this Section we give a short overview of different approaches to
functional calculus.

\subsection{Two Procedures of Functional Calculus}
The following definition is preliminary.
\begin{defn}
A {\em functional calculus\/} $\Phi(\algebra{A},T)$ for an $m$-tuple of
self-adjoint operators $T=(T_1,\ldots,T_m)$ on the Hilbert space $H$
is a continuous mapping $\Phi: \algebra{A} \rightarrow \bos(H)$ of a
linear space \algebra{A} of functions into a linear space of operators on $H$
with
some natural properties.
\end{defn}
We do not specify yet the algebra \algebra{A} and the {\em natural
properties\/}. In many cases the algebra \algebra{A} should at least
contain all smooth functions with compact support (space
$\FSpace{C^\infty}{c}(\Omega)$) in set $\Omega\subset\Space{R}{m} $
depending of the spectra of the operators $T_j$. The natural conditions
usually include some combinations of the following items:
\begin{enumerate}
\item $\Phi(1)=I$, i.~e., the image of the function identically equal to $1$
is the identity operator on $H$;
\item $\Phi(\bar{f})=\Phi(f)^*$;
\item $\Phi(x_j)=T_j$;
\item $\norm{\Phi(f)}_H\,\leq\,\norm{f}_\algebra{A}$;
\item $\Phi(\algebra{A},T)$ depends continuously on the set $T$.
\end{enumerate}
It is notable, that the {\em algebraic\/} properties 1--3 usually imply the
{\em metric\/} ones 4,5.

Generally speaking there are two different but closely connected
ways
to construct a functional calculus.
\begin{procedure}\label{sc:polinom} The homomorphism $\Phi$ is defined on some
dense subalgebra $\algebra{A}_0$ of the
algebra
\algebra{A}. The algebra $\algebra{A}_0$ should have a simple algebraic
structure (for example, the algebra of polynomials), which allows a simple
definition of  $\Phi$. After that $\Phi$ extends to the whole of
\algebra{A} by continuity.
\end{procedure}
\begin{procedure}\label{sc:kernel} Consider a reproducing formula for
functions from \algebra{A} of the form
\begin{equation}
f(x)=\int K(x,y)f(y)\,dy.
\end{equation}
For example, a Cauchy type formula (see
Subsection~\ref{ss:riesz}) or a Fourier type transform (see
Subsection~\ref{ss:weyl}). If the kernel $K(x,y)$ can be extended in a natural
way to an operator $K(T,y)$, then we
 define the functional calculus by the formula
\begin{equation}
\Phi(f)=f(T)=\int K(T,y)f(y)\,dy.
\end{equation}
\end{procedure}
The first Procedure is sometimes easier, but unlike the second one it does
not
give us an explicit formula for functions outside $\algebra{A}_0$.
Usually it is possible to construct, the function calculus in both ways
(see examples bellow). The first Procedure is carried out for the
holomorphic calculus in~\cite{JTaylor72} and the second
one in~\cite{JTaylor70}.

\subsection{Functional calculus: The Weyl approach}\label{ss:weyl}
 Let us recall that a
PDO $\object{Op}a(x,\xi)$~\cite{Hormander85,Shubin87,MTaylor81,HWeyl},
with the {\em
Weyl symbol\/} $a(x,\xi)$, is defined by the formula:
\begin{equation}\label{eq:pdo}
[\object{Op}a](x,\xi)u(y)=\int_{\Space{R}{n}\times\Space{R}{n}}
a(\frac{x+y}{2},\xi)\,e^{i\scalar{y-x}{\xi}}u(x)\,dx\,d\xi.
\end{equation}
This formula may be obtained by Procedure~\ref{sc:kernel} from the
formula
of the inverse
Fourier transform
\begin{displaymath}
f(x)=(2\pi)^{-n/2}\int_\Space{R}{n} \widehat{f}(\xi)e^{ix\xi}\, d\xi.
\end{displaymath}
Namely, let
us take the set of self-adjoint operators
\begin{equation}\label{eq:pdo-frame}
T_j=y_j,\ T_{j+n}=\frac{1}{i}\frac{\partial }{\partial y_j},\ 1\leq
j\leq n,
\end{equation}
which operate on $S=\Space{R}{N}$ in the obvious way. Note that these
operators have exactly the commutators~\eqref{eq:heisen-comm}. Then we
have (see~\cite[\S~1.3]{MTaylor86})
\begin{eqnarray*}
[Kf](y)&=&(2\pi)^{-N}\int_\Space{R}{2N} \widehat{k}(x,\xi)\,
           e^{i(\sum_1^Nx_jy_j-\sum_1^N \xi\frac{\partial }{i\partial
y_j})}f(y)
           \,dx \,d\xi\\
&=&(2\pi)^{-N}\int_\Space{R}{2N} k(\frac{x+y}{2},\xi)\,
           e^{i\scalar{y-x}{\xi}} f(x)\,dx\,d\xi,
\end{eqnarray*}
i.~e.~it defines exactly the {\em Weyl functional calculus\/} (or the
Weyl
quantization). PDO calculus is a very
important tool for the theory of differential equations and quantum
mechanics.
\begin{rem} \label{re:origin}
 Feynman proposed in~\cite{Feynman51} an extension of this functional
calculus---the {\em
functional calculus of ordered operators\/}---in a very similar way.
Anderson~\cite{Anderson69} introduced a generalization of the Weyl
calculus for an arbitrary set $\{T_j\}$ of self-adjoint operators in a
Banach space
by the formula
\begin{equation}
K=(2\pi)^{-N/2}
\int_\Space{R}{N} \widehat{k}(x_1,x_2,\ldots,x_N)\,e^{i\sum_1^N x_j
T_j}\,dx,
\end{equation}
 A description of several different
operator calculus may be found in~\cite{Maslov73}. It was shown by
R.~Howe~\cite{Howe80a,Howe80b} that the success of the original Weyl
calculus is intimately connected with the structure of the Heisenberg
group and its different representations.
\end{rem}
\begin{rem}\label{re:weyl-polin} It should be noted, that besides the
definition of the Weyl functional calculus in accordance with the
Procedure~\ref{sc:kernel} it can be defined also by the
Procedure~\ref{sc:polinom}.
As subalgebra $\algebra{A}_0$ should be taken the algebra of symmetric
 polynomials of the variables $x_j$ and the homomorphism
$\Phi$
``insert'' the operators $T_j$ instead of the
variables $x_j$~\cite[Theorem~2.4]{Anderson69}. We will use this observation
to establish connection between the Weyl and the Riesz-Clifford
calculus.
\end{rem}

\subsection{Functional calculus: The Riesz approach}\label{ss:riesz}
For a pair of self-adjoint operators
$T_1,\ T_2$, the Riesz calculus~\cite[Chap.~XI]{RieszNagy55} can be
also
developed in accordance with the Procedure~\ref{sc:kernel}. The starting
point is the Cauchy integral formula for a function
$f(z)$ holomorphic in $z=x_1+ix_2$ near the spectrum $\sigma(T)$ of the
operator
$T=T_1+iT_2$ and a nice contour $\Gamma$ containing $\sigma(T)$:
\begin{displaymath}
f(x_1+ix_2)=(2\pi)^{-1}\int_\Gamma f(\tau)(\tau -(x_1+ix_2))^{-
1}\,d\tau.
\end{displaymath}
Thus if we are able to define the function $K(z,\tau)=(\tau I-
(x_1+ix_2))^{-
1}$ for an operator $T=T_1+iT_2$, then the functional calculus will be given by
the
formula
\begin{equation}\label{eq:riesz}
f(T_1+iT_2)=(2\pi)^{-1}\int_\Gamma f(\tau)(\tau I-(T_1+iT_2))^{-
1}\,d\tau.
\end{equation}
By Runge's theorem, the Riesz calculus may also be obtained by
Procedure~\ref{sc:polinom}
as the
closure of the polynomials in $T$ in the space of holomorphic functions.
This definition, as was shown
in~\cite[Theorem~5.1]{Anderson69}, gives that the Riesz and the Weyl
calculus are
essentially the
same
in the case of a pair of bounded operators $(T_1,T_2)$.

To generalize the Riesz calculus from a pair of operators to  arbitrary
$m$-tuples, one should use a generalization of complex analysis in one
variable. There are many essentially different theories which can be
considered as generalizations. The principal contributions to
holomorphic functional calculus based on the analysis of several complex
variables were made in~\cite{JTaylor70,JTaylor72}. By the very nature of
complex analysis, this approach can work only for a commuting set
of operators.

Another approach to the Riesz calculus based on Clifford analysis
may be
found in~\cite{McInPryde87}, but there again it was found
 that {\em such a functional calculus is
possible only for an $m$-tuple of commuting generalized scalar
operators
with real spectra\/}~\cite[Theorem~9.1]{McInPryde87}.

Despite of the rich mathematical aspects of these theories, it should be
noted that
quantum mechanics requires a calculus of essentially non-commuting
operators. Moreover, by the von~Neumann
theorem~\cite[\S~90]{AkhGlaz81a} it follows, that for any  $m$-tuple
$\{T_j\}$ of commuting self-adjoint operators, there are such
$m$-tuple
of real valued functions $\{f_j\}$ and an operator $S$ such that
\begin{displaymath}
T_j=f_j(S).
\end{displaymath}
Thus a calculus for $m$-tuple of commuting operators is (in some
sense)
the calculus of only one operator. Therefore one is forced to look
for
another definition of a non-trivial functional calculus of a set of
non-commuting  operators.

\section{Riesz-Clifford Calculus for Non-Commuting
Operators}\label{se:riesz}
To extend the Riesz calculus
to an arbitrary $m$-tuples of bounded operators $\{T_j\}$ it seems
natural
to use Clifford analysis (see for example~\cite{DelSomSou92}), which
is a non-commutative analogy to one-dimensional complex analysis. Then
one can define a
function from an arbitrary $m$-tuple of bounded self-adjoint operators
$\{T_j\}$, either by the Cauchy integral formula (compare
with~\eqref{eq:riesz} and~\cite[Remark~4.11]{Kisil94e}) or by the subset of
symmetric
hyperholomorphic
polynomials~\cite{Malonek93}.

In this Section we will systematically develop such a point of view.

\subsection{Clifford Algebras and Clifford Analysis}
We need some notation from~\cite{DelSomSou92}. Let the euclidean vector
space \Space{R}{n+1} have the orthonormal basis $e_0,e_1,\dots,e_n$. The
Clifford
algebra \Cliff{n}{0} is generated by the elements $e_0,e_1,\dots,e_n$ with
the usual vector operations and the multiplication defined on the elements of
the basis by the
following equalities
\begin{eqnarray*}
e_j^2 &=& -e_0 \qquad (j=1,\ldots,n),\\
e_je_k+e_ke_j&=&0 \qquad (j,k=1,\ldots,n; j\neq k),
\end{eqnarray*}
and then extended linearly to the whole space.
An element of \Cliff{n}{0} can be written as a linear
combination with coefficients $a_\alpha \in \Space{R}{} $ of the
monomials
$e_\alpha = e_{1}^{j_{1}}e_{2}^{j_{2}} \cdots e_{m}^{j_{m}}$ :
\begin{equation}
a=\sum_\alpha  a_\alpha  e_\alpha = \sum_{j_{k}= 0\ or\ 1}
a_{j_{1} j_{2} \ldots j_{m}}
e_{1}^{j_{1}}e_{2}^{j_{2}} \cdots e_{m}^{j_{m}}.
\label{eq:elem}
\end{equation}
The conjugate  $\bar{a}$ of an element $a$ is defined by the rule:
\begin{equation}
\bar{a}=\sum_\alpha  a_\alpha  \bar{e}_\alpha = \sum_{j_{k}= 0\ or\ 1}
a_{j_{1} j_{2} \ldots j_{m}}
\bar{e}_{m}^{j_{m}}\bar{e}_{m-1}^{j_{m-1}} \cdots \bar{e}_{1}^{j_{1}},
\label{eq:conjug-elem}
\end{equation}
where $\bar{e}_j=-e_j,\ \bar{e}_0=e_0,\ 1\leq j\leq n$. A Clifford
algebra valued function $f$
 of the variables $(x_0,x_1, \dots, x_n)$, is called
{\em hyperholomorphic\/} in an open domain $\Omega\subset\Space{R}{n+1} $ if it
satisfies the
Dirac
(or Weyl~\cite{DelSomSou92})  equation
\begin{equation}\label{eq:dirac}
Df=\sum_{j=0}^{n}e_j\frac{\partial f}{\partial x_j}=0.
\end{equation}
The main results of Clifford analysis (Cauchy theorem,
Cauchy integral formula etc.) have a structure closer that of
the complex analysis of one variable than to standard complex analysis of
several variables.

Note one very important feature. In Clifford
analysis {\em not all polynomials are hyperholomorphic functions\/}. Instead
one has to consider the symmetric polynomials of
the monomials having the form
\begin{equation}\label{eq:monom}
\vec{x}_j=(e_jx_0-e_0x_j),\ 1\leq j\leq n.
\end{equation}
The role of such monomials (``regular variables''~\cite{Delanghe70})
is described for quaternionic analysis
in~\cite{Sudbery79}, for Clifford analysis in~\cite{Malonek93},
for Fueter-Hurwitz analysis in~\cite{KrolRam92}, and for solutions of
the general Dirac type equation in~\cite{Kisil94e,Laville87}.

We
use the following facts and notation. Let $\FSpace{H}{}(\Omega)$  denote
 the space of all hyperholomorphic functions in the domain $\Omega$ and
by $\FSpace{P}{}$ the space of all hyperholomorphic polynomials. The
space $\FSpace{P}{}$ has the linear subspaces $\FSpace{P}{j},\ 0\leq j <
\infty$ of homogeneous polynomials degree $j$. We will
show
that $\FSpace{P}{}$ consist of symmetric polynomials
constructed from the
monomials of the form~\eqref{eq:monom} by symmetric
products~{\cite{Malonek93}}
\begin{equation}\label{eq:permut}
a_1\times a_2\times \cdots \times a_k=\frac{1}{k!}\sum a_{j_1} a_{j_2}
\cdots  a_{j_k},
\end{equation}
where the sum is taken over all of permutations of $({j_1} ,{j_2}, \ldots,
{j_k})$. Clifford valued coefficients are written on the right-hand
side.
\begin{lem}
The linear subspace $\FSpace{P}{j}$ has a basis consisting of symmetric
polynomials.
\end{lem}
\begin{proof}
Let $V_\alpha(x) $ be the following homogeneous hyperholomorphic polynomial of
degree $\modulus{\alpha }=k$ (see~\cite[Chap.~II,
Definition~1.5.1]{DelSomSou92}):
\begin{equation}\label{eq:v-polin}
V_\alpha (x)=\sum_{j=0}^{\modulus{\alpha}}\,\frac{(-1)^j
x_0^j}{j!}\,[(\sum_{l=1}^n \frac{\partial }{\partial
x_l})^j\,(\sum_{l=1}^n
x_l)^\alpha].
\end{equation}
Then, the $V_\alpha$ with $\modulus{\alpha}=k $, form a basis of $V_k$. It is
easy to check that
\begin{equation}\label{eq:v-permut}
V_\alpha (x)=(-1)^{\modulus{\alpha}}(e_0x_1-e_1x_0)^{\alpha_1}\times(e_0x_2-
e_2x_0)^{\alpha_2}\times\cdots\times(e_0x_n-e_n
x_0)^{\alpha_n}=\vec{x}^\alpha.
\end{equation}
This gives the assertion.
\end{proof}

Let
\begin{equation}\label{eq:cauchy-ker}
E(y-
x)=\frac{\Gamma(\frac{n+1}{2})}{2\pi^{(m+1)/2}}\,\frac{\overline{y-
x}}{\modulus{y-x}^{n+1}}
\end{equation}
be the Cauchy kernel~\cite[p.~146]{DelSomSou92}  and
\begin{equation}
d\sigma=\sum_{j=0}^n (-1)^j e_j dx_0 \wedge \ldots \wedge[dx_j] \wedge
\ldots \wedge dx_m.
\end{equation}
be the differential form of the ``oriented surface
element''~\cite[p.~144]{DelSomSou92}. Then for any
$f(x)\in\FSpace{H}{}(\Omega)$ we have the Cauchy integral
formula~\cite[p.~147]{DelSomSou92}
\begin{equation}
\int_{\partial \Omega} E(y-x)\,d\sigma_y\,
f(y)=\left\{\begin{array}{l}
f(x),\ x\in\Omega\\
0, x\not\in\bar{\Omega}
\end{array}.\right.
\end{equation}
We should point out the universality (with respect to domains) of both
the
Cauchy kernel $E(y-x)$ and the Cauchy formula in the Clifford analysis,
in contrast to the case of several complex variables.

If we define as in~\cite[\S~18.6]{BraDelSom82} and~\cite[Chap.~II,
Definition~1.5.5]{DelSomSou92}
\begin{equation}\label{eq:w-def}
W^{(a)}_\alpha (x)=(-1)^{\modulus{\alpha}}\partial ^\alpha E(x-a),
\end{equation}
then, for $\modulus{x} < \modulus{y} $ we obtain~\cite[Chap.~II,
(1.16)]{DelSomSou92}
\begin{equation}\label{eq:decomp}
E(y-x)=\sum_{j=0}^\infty\left(\sum_{\modulus{\alpha }=j}V_\alpha
(x)W_\alpha (y)\right).
\end{equation}
Let $M_{(r)}(\Omega)$ be the unitary right Clifford module of all left
hyperholomorphic functions in $\Omega$. For fixed $a$ $\FSpace{R}{(r)}(a) $
denotes the right Clifford module generated by functions $W^{(a)}_\alpha (x)$.
If $S=\{a_i \such i\in \Space{N}{} \}$ is a finite or countable subset of
\Space{R}{n+1}, put $\FSpace{R}{(r)}(S)=\cup_{i\geq 1} \FSpace{R}{(r)}(a_i)$.

The following result plays an important role latter.
\begin{thm}[Runge type] \textup{\cite[\S~18]{BraDelSom82}}\label{th:runge} Let
$K$ be a compact subset of \Space{R}{n+1} and let $S$ be a subset of
$\Space{R}{n+1}\setminus K$ having one point in each bounded component of
$\Space{R}{n+1}\setminus K$. The
n each function which is hyperholomorphic in a neighborhood  of $K$ can be
approximated uniformly on $K$ by the functions in $\FSpace{M}{(r)}
(\Space{R}{n+1})\oplus \FSpace{R}{(r)}(S)$, i.~e.  $\FSpace{M}{(r)}
(\Space{R}{n+1})\oplus \FSpace{R}{(r)}(S)$ is
 uniformly dense in $\FSpace{H}{}(K)$.
\end{thm}

\subsection{Riesz-Clifford calculus: a Polynomial Procedure}

Fix an $m$-tuple of bounded self-adjoint operators
$T=(T_1,\ldots,T_m)$ on the Hilbert space $H$. It is well
known~\cite[Chap.~I]{DelSomSou92} that any Clifford algebra
(particularly, \Cliff{n}{0}) can be realized as an algebra of
endomorphisms of some finite-dimensional vector space $H_0$. Let
$\widetilde{H}=H\otimes H_0$. Then if $n\geq m$ we can associate with $T$  an
operator $\widetilde{T}:\widetilde{H}\rightarrow\widetilde{H}$
by the formula
\begin{displaymath}
\widetilde{T}=\sum_{j=1}^m T_j\otimes e_j.
\end{displaymath}
We use the notation $\bos(H)$ and $\bos(\widetilde{H})$ for the algebra of
bounded operators in the corresponding Hilbert space.
\begin{defn}
Let \algebra{A} be an algebra with the operations of addition
$+_\algebra{A}$ and multiplication $\cdot_\algebra{A}$ which is
generated by a finite set of elements $a_1,\ldots,a_k$.
Define a new operation $\times_\algebra{A}$ of symmetric
multiplication associated with them by the formula~\eqref{eq:permut}.
The resulting set will be called an
{\em $\times$-algebra\/}. Let \algebra{A} and \algebra{B} be two
$\times$-algebras. We say that $\phi:
\algebra{A} \rightarrow \algebra{B}$ is an {\em $\times$-homomorphism\/}
of two $\times$-algebras if  the following holds
\begin{enumerate}
\item $\phi(a_j)=b_j,\ 1\leq j\leq k$, where $a_1,\ldots,a_k$ and
$b_1,\ldots,b_k$ are generators of \algebra{A} and \algebra{B} correspondingly.
\item $\phi(\lambda a_1+a_2)=\lambda\phi(a_1)+\phi(a_2)$ for any $a_1,
a_2\in\algebra{A}$ and $\lambda\in \Space{C}{}$.
\item $\phi(a_1\times\cdots\times a_n)=\phi(a_1)\times\cdots\times
\phi(a_n)$ for any set $a_1,\ldots,a_n$ of (not necessarily distinct)
generators of \algebra{A}.
\end{enumerate}
\end{defn}
By the definition of the symmetric product one can obviously
deduce
\begin{lem}\label{le:homom} Any homomorphism of two algebras is a
$\times$-homomorphism of the corresponding $\times$-algebras.
\end{lem}

The converse, of course, is not true:
\begin{example}\label{ex:homom}
Let \algebra{A} be the algebra of polynomials in two variables
$x_1$ and $x_2$, and \algebra{B} be the algebra generated by two
operators $P$ and $Q$ with the commutator $[P,Q]=i\hbar I$. Then the mapping
$\phi$ defined on the generators of algebra \algebra{A} by
\begin{displaymath}
 \phi:x_1\mapsto P,\ \phi: x_2 \mapsto Q
\end{displaymath}
cannot be extended to an homomorphism of algebras \algebra{A} and
\algebra{B} but defines a $\times$-homomorphism, which coincides with the
Weyl quantization.
\end{example}

The following is a definition of Riesz-Clifford hyperholomorphic calculus
(compare with~\cite{McInPryde87}).
\begin{defn}\label{de:calculus} We say that
$\widetilde{T}$ {\em has a functional calculus $(\algebra{A}, \Phi)$ based
on
\Space{R}{m}\/} whenever the following conditions hold: \algebra{A} is a
topological vector space of hyperholomorphic functions from
\Space{R}{m} to $\Cliff{0}{n}$ ($m\leq n$), with addition
and (symmetric) multiplication defined pointwise, and $\Phi:
\algebra{A}\rightarrow
\bos(\widetilde{H})$ is a continuous mapping and its restriction to
symmetric polynomials is a $\times$-homomorphism such that
\begin{equation}\label{eq:phi-def}
\Phi: \vec{x}_j (=e_j x_0 - e_0 x_j) \mapsto T_j,\  1\leq j \leq m
\end{equation}
\end{defn}
There, to extend calculus from commuting operators to
non-commuting we ones we  relax the requirement from  homomorphism
to  $\times$-homomorphism.
\begin{defn}\label{de:support}
The {\em support\/} of the homomorphism $\Phi$ is the smallest closed set in
\Space{R}{n} such that $\Phi(f)=0$ for all $f\in \algebra{A}$ with
$\spt\Phi\cap\spt f=\emptyset$. The {\em joint $\times$-resolvent set\/}
$R_\times(T)$ of $T$ is the largest open subset of \Space{R}{m} such
that for any point $\lambda=(\lambda_1,\ldots,\lambda_m)\in
R_\times(T)$, the $\times$-algebra generated by the operators
$T_j\otimes e_j-\lambda_j I\otimes e_j, 1\leq j\leq m$ contains at least
one invertible element in $\widetilde{H}$. The {\em joint spectrum\/}
$\sigma_\times(T)$ of the operators $T_1,
\ldots,T_m$ is the set $\Space{R}{m}\setminus \sigma_\times(T)$.
\end{defn}
\begin{thm}[Uniqueness]
In the case of $n=m$ for a given simply-connected domain $\Omega$ and an
$m$-tuple of operators $T$, there exists no more than one hyperholomorphic
functional calculus.
\end{thm}
\begin{proof}
By the $\times$-homomorphism conditions one can extend the
correspondence~\eqref{eq:phi-def} to the space $\FSpace{P}{}$ of
hyperholomorphic polynomials in only one way. Then the denseness of
$\FSpace{P}{}$ in $\FSpace{H}{}(\Omega)$ implies the assertion.
\end{proof}
\begin{thm}
For any $m$-tuple $T$ of bounded self-adjoint operators there exist a
hyperholomorphic calculus on \Space{R}{m}.
\end{thm}
\begin{proof}
Any hyperholomorphic function on $f(x)\in \FSpace{H}{}(\Space{R}{m})$ can be
represented as the Taylor series~\cite[Chap.~II, \S~1.6.3,
Remark~2]{DelSomSou92}
\begin{displaymath}
f(x)=\sum_{\modulus{\alpha }=1}^\infty V_\alpha (x) c_\alpha,\ c_\alpha \in
\Cliff{n}{0}.
\end{displaymath}
We have already seen what is $V_\alpha (T)$, thus the functional calculus may
be defined by the formula
\begin{displaymath}
f(T)=\sum_{\modulus{\alpha }=1}^\infty V_\alpha (T) c_\alpha.
\end{displaymath}
The convergence of this series is evident.
\end{proof}

\subsection{Riesz-Clifford Calculus: an Integral Representation
Procedure}

We now deduce the explicit formula for hyperholomorphic calculus. Accordingly
to Procedure~\ref{sc:kernel} we should define the Cauchy
kernel
of the operators $T_j$ as follows.
\begin{defn}
Let $n\geq m$. Then (cf.~\ref{eq:decomp})
\begin{equation}\label{eq:cauchy-op}
E(y,T)=\sum_{j=0}^\infty\left(\sum_{\modulus{\alpha }=j}V_\alpha
(T)W_\alpha (y)\right)
\end{equation}
where
\begin{equation}
V_\alpha(T) = T_1^{\alpha _1}\times\cdots \times T_m^{\alpha _m}\times
I\times\cdots\times I;
\end{equation}
i.~e., we have formally substituted in~\eqref{eq:v-permut} for
$(\vec{x}_1,\ldots,\vec{x}_n)$ the $n$-tuple of operators
$(T_1,\ldots,T_m,I,\ldots,I)$.
\end{defn}
\begin{rem}
For the sake of the simplicity we will consider in this Section only
the
case  $m=n$. Consideration of the  more general case $n\geq m$ is
analogous and will be needed Section~\ref{se:field}.
\end{rem}
\begin{lem}\label{le:bound}
Let $\modulus{T}=\max_j\{\norm{T_j}\}$. Then for fixed
$\modulus{y}>\modulus{T}$, the equation~\eqref{eq:cauchy-op} defines a
bounded operator in $\widetilde{H}$.
\end{lem}
\begin{proof}
Consideration of the real valued series
\begin{eqnarray*}
\norm{E(y,T)}&\leq&\sum_{j=0}^\infty\left(\sum_{\modulus{\alpha
}=j}\norm{V_\alpha (T)}\,\norm{W_\alpha (y)}\right)\\
&\leq&\sum_{j=0}^\infty\left(\sum_{\modulus{\alpha
}=j}\modulus{T}^j\,\norm{W_\alpha (y)}\right),
\end{eqnarray*}
and the properties of the Cauchy kernels shows that for the mentioned $y$
the
series~\eqref{eq:cauchy-op} converges in the uniform operator topology to a
bounded operator.
\end{proof}
\begin{defn}\label{de:spectra}
The maximal open subset $R_C(T)$ of $\Space{R}{n}$ such that for any
$y\in R_C(T)$ the series in~\eqref{eq:cauchy-op}
converges in the uniform operator topology to a bounded operator on
$\widetilde{H}$ will be called the {\em (Clifford) resolvent set\/} of $T$.
The complement of $R_C(T)$ in \Space{R}{n} will be called the {\em (Clifford)
spectral set\/} of $T$ and denoted by $\sigma_C(T)$.
\end{defn}
{}From Lemma~\ref{le:bound} it follows that $R_C(T)$ is always non-
empty
and from Definition~\ref{de:spectra} one can see that
$\sigma_C(T)$
is closed. Moreover, it is easy to see that
\begin{lem}
The Clifford spectral set $\sigma_C(T)$ is compact.
\end{lem}
\begin{proof}
$\sigma_C(T)$ is a closed subset of \Space{R}{n} which is bounded due to
Lemma~\ref{le:bound}.
\end{proof}
\begin{lem}\label{le:v-lem}
Let $r>\modulus{T}$ and let $\Omega$ be the ball $\Space{B}{}(0,r)\in
\Space{R}{m}$. Then for any symmetric polynomial $P(\vec{x})$ we
have
\begin{equation}\label{eq:int-polin}
P(T)=\int_{\partial \Omega} E(y,T)\, d\sigma_y\, P(y)
\end{equation}
where $P(T)$ is the symmetric polynomial of the $m$-tuple $T$.
\end{lem}
\begin{proof}
First, consider the case of the polynomial
$P(\vec{x})=V_\alpha
(\vec{x})$. Recall the formula~\cite[Chap.~II, Lemma
1.5.7(i)]{DelSomSou92}
\begin{equation}
\int_{\partial \Space{B}{}(0,r)} W_\beta(y)\, d\sigma\, V_\alpha
(y)=\delta _{\alpha \beta}.
\end{equation}
Let  $A$ be an operator defined by the right-hand side
of~\eqref{eq:int-polin}, since the convergence in~\eqref{eq:cauchy-op}  is
uniform,
we obtain
\begin{eqnarray*}
A&=&\int_{\partial \Omega} E(y,T)\, d\sigma_y\, P(y)\\
  &=&\int_{\partial \Omega} E(y,T)\, d\sigma_y\, V_\alpha (y)\\
  &=&\int_{\partial \Omega}
  \sum_{j=0}^\infty\left(\sum_{\modulus{\beta}=j}V_\beta(T)
      W_\beta(y)\right) \, d\sigma_y\, V_\alpha (y)\\
  &=& \sum_{j=0}^\infty\left(\sum_{\modulus{\beta }=j}V_\beta (T)\,
  \int_{\partial \Omega} W_\beta (y)\, d\sigma_y\, V_\alpha
(y)\right)\\
  &=& \sum_{j=0}^\infty\left(\sum_{\modulus{\beta }=j}V_\beta
  (T) \delta _{\alpha \beta}\right)\\
  &=& V_\alpha (T)\\
  &=& P(T).
\end{eqnarray*}
Finally, by linear extension we get the assertion.
\end{proof}

\begin{lem}
For any domain $\Omega$ which does not contain $\sigma_C(T)$ and any
$f\in\FSpace{H}{}(\Omega)$, we have
\begin{equation}
\int_{\partial \Omega} E(y,T)\, d\sigma_y\, f(y)=0
\end{equation}
\end{lem}
\begin{proof}
If we apply a bounded linear functional $w$ on the space
$\bos(\widetilde{H})$ to the integral at left-hand side of the assertion,
we will obtain an integral of the hyperholomorphic function
\begin{displaymath}
w(E(y,T))\, d\sigma_y\, f(y),
\end{displaymath}
which is equal to $0$ by the Cauchy theorem. The arbitrariness of $w$
gives the assertion.
\end{proof}

Due to this Lemma we can replace the domain  $\Space{B}{}(0,r) $ at
Lemma~\ref{le:v-lem} with an arbitrary domain $\Omega$
containing the spectral set $\sigma_C(T).$
An application of Lemma~\ref{le:v-lem} gives the
main
\begin{thm}
Let $T=(T_1,\ldots,T_m)$ be an $m$-tuple of bounded self-adjoint
operators. Let the domain $\Omega$ with piecewise smooth boundary have a
connected complement and suppose the spectral set $\sigma_C(T)$ lies inside a
domain $\Omega$. Then the mapping
\begin{equation}\label{eq:int-calc}
\Phi: f(x)\mapsto f(T)=\int_{\partial \Omega} E(y,T)\, d\sigma_y\,
f(y)
\end{equation}
defines a hyperholomorphic calculus for $T$.
\end{thm}
\begin{proof}
We have already seen (Lemma~~\ref{le:v-lem}) that formula~\eqref{eq:int-calc}
defines a
functional calculus for functions in $\FSpace{M}{(r)}
(\Space{R}{n+1})$. Application of
Theorem~\ref{th:runge} allows us to extend this calculus to the whole of
$\FSpace{H}{}(\Omega)$.
\end{proof}

\section{Riesz-Clifford Calculus and Quantum Field
Theory}\label{se:field}
\subsection{Quantum Field Theory and Clifford
Algebras}\label{ss:qfield}

The step from a  quantum system with one particle to quantum ensemble
with $m$ particles depends on the kind of particles we have. For
example, photons satisfy to Bose-Einstein
statistics, but electrons satisfy to Fermi-Dirac
statistics \cite{Dirac67}. This means that observables of a system of photons
are
symmetrical operators relative to equivalent particles; but
observables
of a system of electrons are anti-symmetric ones. This involves, in
particular, the Pauli's exclusion principle: there is only one or no
Fermi particle in each state. Such a distribution of
particles among different states can be easily described with the help
of the Clifford algebra \Cliff{n}{0}. The monomial
$e_{1}^{i_{1}}e_{2}^{i_{2}} \cdots e_{m}^{i_{m}}$
can be considered as the description of a physical system, in
which the only states enumerated by $k$, for which $i_{k}\neq
0$, are filled by a Fermi particle. The general element
\begin{displaymath}
a=\sum_\alpha  a_\alpha  e_\alpha = \sum_{j_{k}= 0\ or\ 1}
a_{j_{1} j_{2} \ldots j_{m}}
e_{1}^{j_{1}}e_{2}^{j_{2}} \cdots e_{m}^{j_{m}}
\end{displaymath}
 of \Cliff{m}{0} can be interpreted as
a probability distribution among  ``pure'' states defined by such
monomials\footnote{We consider only the case of a finite number
of different states. The difficulties occurred in the infinite case is
described for example in~\cite{Segal94}.}~\cite{Dirac67}.

\subsection{Quantum Field Theory and Riesz-Clifford Calculus}
Non-commutativity of quantum observables generate the following
question: {\em Let $A$ and $B$ be quantum analogues of classic
observables $a$ and $b$. What corresponds to the classic observable
$ab$, the result of simultaneous measurement of $a$ and $b$?\/}

Jordan (see for example~\cite[\S~1.2]{Kirillov90}) proposed the answer
\begin{equation}\label{eq:jordan}
A\circ B= \frac{1}{2}(AB+BA)=
\left(\frac{A+B}{2}\right)^2-\left(\frac{A-B}{2}\right)^2.
\end{equation}
One can interpret this formula in the following way. Simultaneous
measurement of $A$ and $B$ from the {\em macro\/}scopic point of view is not
simultaneous from the {\em micro\/}scopic one. It is natural to assume that two
different processes (measurement of $AB$ and $BA$) will occur with equal
probability $\frac{1}{2}$
. This is explicitly described\footnote{ In spite of some
advantages of the Jordan algebras, they were not too successful and
have almost disappeared in the modern literature.}  by~\eqref{eq:jordan}.

To prolong such a point of view for three operators $A$, $B$, $C$ we can accept
neither $(A\circ B)\circ C $ nor $A\circ (B\circ C)$ (the Jordan
multiplication~\eqref{eq:jordan} is not associative!). We again should select
the equal probability $\frac{1}{
6}$ for six different processes $ABC$, $ACB$, $BAC$, $BCA$, $CAB$, $CBA$ and
take the value
\begin{displaymath}
\frac{1}{6}(ABC+ACB+BAC+BCA+CAB+CBA)=A\times B\times C.
\end{displaymath}
Note also that $A\circ B=A\times B$.

Prolonging this line we should agree that the symmetric monomial
\begin{displaymath}
A_1\times A_2\times \cdots \times A_k
\end{displaymath}
correspond to the monomial of classic observables
\begin{displaymath}
a_1\times a_2\times \cdots \times a_k.
\end{displaymath}
Thus the symmetric polynomials are a quantum version of  classic polynomials.
Then the natural condition of continuity leads us to the Riesz-Clifford
functional calculus.

One can summarize the main point of this Section the following
way: To describe of an ensemble of Fermi particles we have to study
the algebra of operators on the space $\widetilde{H}=H\otimes H_0$,
where the Hilbert space $H$ is a space of states of one particle and
$H_0$ is a (finite-dimensional in our consideration) space of filling
numbers. The Riesz-Clifford functional calculus based on the symmetric
product naturally describes functions of non-commuting quantum observables in
quantum field theory.

\makeatletter \renewcommand{\@biblabel}[1]{\hfill#1.}\makeatother
\newcommand{\bysame}{\leavevmode\hbox to3em{\hrulefill}\,}
\newcommand{\noopsort}[1]{} \newcommand{\printfirst}[2]{#1}
  \newcommand{\singleletter}[1]{#1} \newcommand{\switchargs}[2]{#2#1}
  \newcommand{\irm}{\mbox{\rm I}} \newcommand{\iirm}{\mbox{\rm II}}
  \newcommand{\vrm}{\mbox{\rm V}}


\begin{thebibliography}{10}

\bibitem{AkhGlaz81a}
N.~I. Akhiezer and I.~M. Glazman.
\newblock {\em Theory of Linear Operators in {Hilbert} Space}, volume~1.
\newblock Pitman Advanced Publishing Program, London, 1981.

\bibitem{Anderson69}
R.~F.~V. Anderson.
\newblock The {Weyl} functional calculus.
\newblock {\em J.~Funct. Anal.}, 4:240--267, 1969.

\bibitem{Berezin74}
F.~A. Berezin.
\newblock Quantization.
\newblock {\em Math. USSR Izvestija}, 8:1109--1165, 1974.

\bibitem{BraDelSom82}
F.~Brackx, R.~Delanghe, and F.~Sommen.
\newblock {\em Clifford Analysis}, volume~76 of {\em Research Notes in
  Mathematics}.
\newblock Pitman Advanced Publishing Program, Boston, 1982.

\bibitem{Coburn94a}
L.~A. Coburn.
\newblock {Berezin-Toeplitz} quantization.
\newblock In {\em Algebraic Mettods in Operator Theory}, pages 101--108.
  Birkh\"auser, New York, 1994.

\bibitem{CobXia94}
L.~A. Coburn and J.~Xia.
\newblock Toeplitz algebras and {Rieffel} deformation.
\newblock {\em Comm. Math. Phys.}, 1994.
\newblock to appear.

\bibitem{Delanghe70}
R. Delanghe.
\newblock On the singularities of functions with values in a
           Clifford algebra.
\newblock {\em Math. Ann.\/} {\bf 196}(1972), pages 293--319.


\bibitem{DelSomSou92}
R.~Delanghe, F.~Sommen, and V.~Sou\v{c}ek.
\newblock {\em Clifford Algebra and Spinor-Valued Functions}.
\newblock Kluwer Academic Publishers, Dordrecht, {\noopsort{}}1992.

\bibitem{Dirac67}
P.~A.~M. Dirac.
\newblock {\em Lectures on Quantum Field Theory}.
\newblock Yeshiva University, New York, {\noopsort{}}1967.

\bibitem{CAlgebras94}
R.~Doran, editor.
\newblock {\em {$C^*$}-Algebras: 1943--1993}.
\newblock Number 167 in Contemporary Mathematics. AMS, Providence, Rhode
  Island, 1994.

\bibitem{Effros94}
E.~G. Effros.
\newblock Some quantizations and reflections inspired by the {Gelfand-Naimark}
  theorem.
\newblock In Doran \cite{CAlgebras94}, pages 99--113.

\bibitem{Feynman51}
R.~P. Feynman.
\newblock An operator calculus having applications in quantum electrodynamics.
\newblock {\em Phys.~Rev.}, 84(2):108--128, 1951.

\bibitem{Folland89}
G.~B. Folland.
\newblock {\em Harmonic Analysis in Phase Space}.
\newblock Princeton University Press, Princeton, New Jersey, {\noopsort{}}1989.

\bibitem{Hormander85}
L.~H\"{o}rmander.
\newblock {\em The Analysis of Linear Partial Differential Operators
  \irm\irm\irm: {P}seudodifferential Operators}.
\newblock Springer--Verlag, Berlin Heidelberg New York Tokyo, 1985.

\bibitem{Howe80a}
R.~Howe.
\newblock On the role of the {Heisenberg} group in harmonic analysis.
\newblock {\em Bull. of the AMS (New Series)}, 3(2):821--843, 1980.

\bibitem{Howe80b}
R.~Howe.
\newblock Quantum mechanics and partial differential equations.
\newblock {\em J. Funct. Anal.}, 38:188--254, 1980.

\bibitem{Kirillov90}
A.~A. Kirillov.
\newblock {\em Geometric Quantization}, volume~4 of {\em Encyclopedia of
  Mathematical Sciences}, pages 137--172.
\newblock Springer-Verlag, Berlin, 1990.

\bibitem{Kisil94e}
V.~V. Kisil.
\newblock Relative convolutions. \irm. {Properties} and applications.
\newblock Reporte Interno \#~162, Departamento de Matem\'aticas, CINVESTAV del
  I.P.N., Mexico City, 1994 (to appear in {\em Advances in Mathematics\/}).

\bibitem{KrolRam92}
W.~Kr\'olikowski and E.~Ram\'{\i}rez~de Arellano.
\newblock Polynomial solutions of the {Fueter-Hurwitz} equation.
\newblock In A.~Nagel and E.~L. Stout, editors, {\em The Madison Symposium on
  Complex Analysis}, number 137 in Contemprorary Mathematics, pages 297--305.
  AMS, Providence, Rhode Island, 1992.

\bibitem{Laville87}
G.~Laville.
\newblock Sur un calcul symbolique de {Feynmann}.
\newblock In {\em Seminar d'Analyse}, volume 1295 of {\em Lect. Notes in
  Math.}, pages 132--145. Springer--Verlag, Berlin, {\noopsort{}}1987.

\bibitem{Malonek93}
H.~R. Malonek.
\newblock Hypercomplex differentiability and its applications.
\newblock In F.~et~al., editor, {\em Clifford Algebras and Applications in
  Mathematical Physics}, pages 141--150. Kluwer Academic Publishers,
  Netherlands, {\noopsort{}}1993.

\bibitem{Maslov73}
V.~P. Maslov.
\newblock {\em Operational Methods}.
\newblock ``Nauka'', Moscow, {\noopsort{}}1973.

\bibitem{McInPryde87}
A.~McIntosh and A.~Pryde.
\newblock A functional calculus for several commuting operators.
\newblock {\em Indiana Univ. Math.~J.}, 36:421--439, 1987.

\bibitem{Rieffel94}
M.~A. Rieffel.
\newblock Quantization and {$C^*$}-algebras.
\newblock In Doran \cite{CAlgebras94}, pages 68--97.

\bibitem{RieszNagy55}
F.~Riesz and B.~Sz-Nagy.
\newblock {\em Functional Analysis}.
\newblock Ungar, New York, {\noopsort{}}1955.

\bibitem{Segal94}
I.~Segal.
\newblock {$C^*$}-algebras and quantization.
\newblock In Doran \cite{CAlgebras94}, pages 55--65.

\bibitem{Shubin87}
M.~A. Shubin.
\newblock {\em Pseudodifferential Operators and Spectral Theory}.
\newblock Springer-Verlag, Berlin, 1987.

\bibitem{Sudbery79}
A.~Sudbery.
\newblock Quaternionic analysis.
\newblock {\em Math. Proc. Camb. Phil. Soc.}, 85:197--225, 1979.

\bibitem{JTaylor70}
J.~L. Taylor.
\newblock The analytic-functional calculus for several commuting operators.
\newblock {\em Acta Math.}, 125:1--38, 1970.

\bibitem{JTaylor72}
J.~L. Taylor.
\newblock A general framework for a multioperator functional calculus.
\newblock {\em Adv. Math.}, 9:183--252, 1972.

\bibitem{MTaylor81}
M.~E. Taylor.
\newblock {\em Pseudodifferential Operators}, volume~34 of {\em Princeton
  Mathematical Series}.
\newblock Princeton University Press, Princeton, New Jersey, 1981.

\bibitem{MTaylor86}
M.~E. Taylor.
\newblock {\em Noncommutative Harmonic Analysis}, volume~22 of {\em Math. Surv.
  and Monographs}.
\newblock American Mathematical Society, Providence, Rhode Island,
  {\noopsort{}}1986.

\bibitem{HWeyl}
H.~Weyl.
\newblock {\em The Theory of Groups and Quantum Mechanics}.
\newblock Dover, New York, {\noopsort{}}1950.

\end{thebibliography}
\end{document}